%% file: main.tex
\documentclass[journal]{IEEEtran}

\usepackage{cite}

\ifCLASSINFOpdf
  \usepackage[pdftex]{graphicx}

\else
    \usepackage[dvips]{graphicx}
\fi

\usepackage{amsmath}
\usepackage{algorithmic}
\usepackage{array}

\ifCLASSOPTIONcompsoc
  \usepackage[caption=false,font=normalsize,labelfont=sf,textfont=sf]{subfig}
\else
  \usepackage[caption=false,font=footnotesize]{subfig}
\fi

\ifCLASSOPTIONcaptionsoff
  \usepackage[nomarkers]{endfloat}
 \let\MYoriglatexcaption\caption
 \renewcommand{\caption}[2][\relax]{\MYoriglatexcaption[#2]{#2}}
\fi

\usepackage{url}
\usepackage{cite}
\usepackage{tabularx}
\usepackage{float}
\usepackage{subfig}
\usepackage{comment}

\usepackage{url}
\usepackage{bm}
\usepackage{tikz}
\usepackage{courier}
\usepackage{verbatim}
\usepackage{amssymb}
\usepackage{framed}
\usepackage{listings}
\usepackage{outlines}
\usepackage{glossaries}
\usepackage{kbordermatrix}
\usepackage[binary-units]{siunitx}
\usepackage{mathtools}
\usepackage{algorithmic}
\usepackage{algorithm}
\usepackage{yhmath}
\usepackage{enumitem}
\usepackage[utf8]{inputenc}

\input{commands.tex}

\input{MATLAB_colors.tex}

\begin{document}


\title{An Array Decomposition Method for Finite Arrays with Electrically Connected
Elements for fast Toeplitz Solvers}


\author{Lucas {\AA}kerstedt, 
        Harald Hultin, 
        and~B. L. G. Jonsson
\thanks{Lucas {\AA}kerstedt and B. L. G. Jonsson are with the School of EECS,
  KTH Royal Institute of Technology, 100 44 Stockholm, Sweden (e-mail:
lucasak@kth.se).}
\thanks{Harald Hultin is with Saab Surveillance, 175, 41 J{\"a}rf{\"a}lla,
Sweden, and also with the School of EECS, KTH Royal Institute of Technology, 100 44
Stockholm, Sweden }}


\maketitle

\begin{abstract}

A large part of the geometry of array antennas is often partially defined by
finite translational symmetries. Applying the \gls{MoM} with the RWG-like
element on an appropriately structured mesh to these arrays results in an
impedance matrix where the main part exhibits a multilevel block Toeplitz
structure. This article introduces a memory-efficient construction method that
effectively represents and reuses impedance calculations. The proposed method,
applicable to electrically connected elements, also accounts for all
non-symmetric parts of the array. The core idea involves nine distinct
electrically connectable components from which the array can be assembled. The
derived multilevel block Toeplitz matrix is further utilized by an in-house
inverse solver to achieve faster and more memory-efficient \gls{MoM} current
vector calculations. We demonstrate the method by computing the far-field of a
32x32 array and the scattering parameters of two tightly coupled 9x9 arrays.
This approach reduces the memory allocation from $\mathcal{O}(N_x^2 N_y^2)$ to
$\mathcal{O}(N_x N_y)$, for an $N_x \times N_y$ array.
\end{abstract}


\IEEEpeerreviewmaketitle

\section{Introduction}

\IEEEPARstart{E}{lectrically} large array antennas are becoming increasingly
popular due to their high gain and beam scanning capabilities
\cite{zhang_6g_2023}. It is proposed that these large array antennas, consisting
of over 1000 elements, will form the basis of the next-generation mobile
networks \cite{zhang_6g_2019, saad_vision_2020}. Designing large array antennas
is, however, a challenge as their electromagnetic analysis is computationally
expensive \cite{lucente_iteration-free_2008, michielssen_multilevel_1996}. With
a full-wave analysis tool such as ordinary \gls{MoM}, the memory allocation
scales as $\mathcal{O}(N^2)$, where $N$ is the number of basis functions. This,
in turn, limits ordinary \gls{MoM} for the analysis of electrically large
antennas\cite{helander_comparison_2017}. 

There exist two main categories of methods that are used in conjunction with
\gls{MoM} which reduce the memory allocation: Methods for solving arbitrarily
shaped antennas, and methods specifically using the finite translation symmetry
of array antennas. Another often used approximation is unit-cell analysis
\cite{craeye_efficient_2004,neto_truncated_2000}. In this article, we focus
solely on the full-wave analysis of finite arrays with translation symmetry.

Methods for solving large, arbitrarily shaped antennas include the \gls{FMM}
\cite{coifman_fast_1993}, the \gls{MLFMM} \cite{lu_multilevel_1994}, and the
\gls{DDM} \cite{peng_integral_2011, martin_multiresolution_2024,
solis_accurate_2020}. For the finite element method, \gls{DDM} is widely used
\cite{lu_domain_1997, zhao_domain_2007, liu_fast_2025, lee_non-overlapping_2005,
jiang_efficient_2024}.

Similarly to the partitioning procedure performed in \gls{FMM}, \gls{MLFMM}, and
\gls{DDM}, there exist methods such as the \gls{SFX}
\cite{matekovits_analysis_2007} and the \gls{CBFM}
\cite{prakash_characteristic_2003}. Both of these methods rely on macro basis
functions with support on the subdomains, to reduce the size of the \gls{MoM}
impedance matrix. While the \gls{SFX} and \gls{CBFM} are for solving antennas of
arbitrary shape, they can be modified to efficiently solve antennas with finite
translation symmetry, as is done in \eg,
\cite{maaskant_fast_2008,lu_accurate_2004,xiang_new_2019}.

Iterative methods, exploiting the finite translation symmetry of array antennas,
\cite{bleszynski_block-toeplitz_2003, elizabeth_h, kindt_array_2003,
brandt-moller_extended_2024} utilize that the interaction between the subdomains
yields a multilevel block Toeplitz \gls{MoM} impedance matrix, such that the
necessary memory allocation decreases. Additionally, the Toeplitz
structure allows for fast matrix-vector multiplication using the FFT,
accelerating the iterative inverse calculation\cite{lee_fast_1986}. 

The multilevel block Toeplitz structure arising from the finite translation
symmetry in array antennas has been noted on multiple occasions, \eg, in \cite{
bleszynski_block-toeplitz_2003, kindt_array_2003, maaskant_fast_2008,
helander_synthesis_2018, zhao_hierarchical_2022, brandt-moller_extended_2024},
but seldom for the case of \emph{electrically connected} elements. For arrays
with electrically connected elements, there exist the \gls{HO-ADM}
\cite{brandt-moller_extended_2024}, where the connected elements are handled
using the discontinuous Galerkin method \cite{peng_discontinuous_2013}. There is
currently no (to the best of our knowledge) method for constructing the
multilevel block Toeplitz \gls{MoM} impedance matrix with electrically connected
elements, apart from \cite{brandt-moller_extended_2024}. As the \gls{HO-ADM}
does not employ standard RWG-like basis functions, there is a clear need for a
method for RWG-like basis functions.

In this article, we propose a strategy for fast and memory-efficient
computation of the \gls{MoM} impedance matrix by exploiting the finite
translation symmetry of the array antenna. To handle electrically connected
elements, we represent the array antenna by nine array components: A center
array element component, and eight margin components surrounding it, similar
to \cite{ xiang_new_2019, brandt-moller_extended_2024, jiang_efficient_2024,
lu_accurate_2004}. By applying our proposed partitioning algorithm for
electrically connected elements on the nine-component structure, the impedance
interaction between spatially displaced array components can be calculated. As a
result, the majority of the \gls{MoM} impedance matrix inherits a multilevel
block Toeplitz structure. Additionally, we show that the interaction between
the margin components and the element components yields block Toeplitz
structured matrices, apparent in the outer regions of the \gls{MoM} impedance
matrix. This methodology reduces the necessary memory allocation by an order of
magnitude, enabling the calculation of large arrays with electrically connected
array elements, \eg, arrays with elements over a finite ground plane. We also
show how the nine-component representation can be used to calculate
the far-field effectively. Furthermore, we demonstrate that with the nine component
array representation, multiple same-sized arrays with non-identical array
elements can be modeled, such that the mutual coupling between two large arrays
can be calculated.

This article is organized as follows: The theory of the proposed array
decomposition method is discussed in Section \ref{sec:Theory}. In Section
\ref{sec:Numerical}, the proposed method is used to calculate the scattering
parameters and the far-field of various large arrays and their memory scaling.
Section \ref{sec:Conclusions} concludes the article.

\section{Theory}
\label{sec:Theory}
The work presented in this article is based on the \gls{EFIE}
\cite{rao_electromagnetic_1982} applied to surface currents $\bm{J}(\bm{r})$ on
a surface $S$, utilizing the RWG-basis $\bm{f}(\bm{r})$:
\begin{equation}
  \bm{J}(\bm{r}) \approx \sum_{n = 1} ^{N} I_n \bm{f}_n(\bm{r}) \text{,}
  \label{eq:currentSum}
\end{equation}
where, $I_n\in\mathbb{C}$ is the current coefficients for the basis functions
$\bm{f}_n(\bm{r})$, which is associated with an internal edge $e_n$, shared by two
supporting triangles $T_n^+$ and $T_n^-$ such that $T_n =  T_n^+ + T_n^-$ is the
support of the $\bm{f}_n(\bm{r})$.

The matrix representation of the \gls{EFIE} is 
\begin{equation}
  Z I = V \text{,}
  \label{eq:MoMMatrixEq}
\end{equation}
where $V$ is the voltage vector, $V\in\mathbb{C}^{N}$. Each element in the
\gls{MoM} impedance matrix $Z$ is defined by
\begin{equation}
  \begin{aligned}
    z_{mn} =& \\
    \frac{-\mathrm{j} \eta}{k} & \int_{T_n}
    \int_{T_m} \nabla_S' \cdot \bm{f}_n(\bm{r}') G(\bm{r},\bm{r}')
      \nabla_S \cdot \bm{f}_m(\bm{r}) \mathrm{d}S \mathrm{d}S' \\
			       & + \mathrm{j} k \eta 
      \int_{T_n}
      \int_{T_m} G(\bm{r},\bm{r}')
      \bm{f}_n(\bm{r}')  \cdot \bm{f}_m(\bm{r}) \mathrm{d}S \mathrm{d}S' \text{,}
  \end{aligned}
  \label{eq:MoMimpedanceMatrix}
\end{equation}
where $k$ is the wavenumber, $\eta$ is the free space impedance. The tangential
divergence operator, $\nabla_S$, is defined by $\nabla = \hat{\bm{n}} \partial_n
+ \nabla_S$, where $\hat{\bm{n}}$ is the normal to $S$. The standard Green’s
function $G(\bm{r},\bm{r}')$, in \eqref{eq:MoMimpedanceMatrix} is defined as $
G(\bm{r},\bm{r}') = \exp{(-\mathrm{j}k|\bm{r} - \bm{r}'|)}/{4\pi |\bm{r} -
\bm{r}'|}$.

\subsection{Subdomain Partitioning in \glsentryshort{MoM}}
The \gls{MoM} impedance matrix \eqref{eq:MoMimpedanceMatrix} is symmetric, such
that it suffices to calculate the upper triangular part. Similarly, for an
impedance matrix partitioned into $P\times P$ submatrices, it suffices to
calculate the upper triangular of blocks.

The submatrices (or blocks) are obtained by partitioning the edge elements into
$P$ subdomains. A submatrix $Z_{ij}$ then corresponds to the impedance relation
between edges of subdomain $i$ and $j$. 
\begin{figure}[ht!]
  \begin{center}
    \includegraphics[width=\columnwidth]{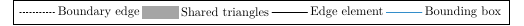}

    \subfloat[]{\includegraphics[width=100pt]{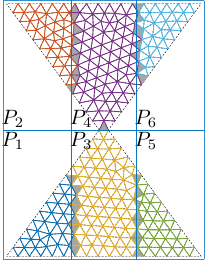}}
    \hfill
    \subfloat[]{\includegraphics[width=100pt]{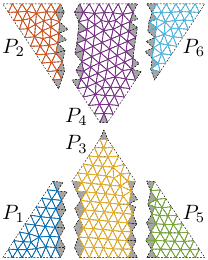}}
  \end{center}
  \caption{Bowtie mesh partitioned into six subdomains with bounding boxes
  $P_1,\dots,P_6$. (a) The partitioned mesh and the bounding boxes. (b) The
  exploded view, where shared triangles are highlighted.}
  \label{fig:meshPartition1}
\end{figure}

In this work, we partition the edge elements according to their physical
position. Edges whose geometric centers are located within a bounding box $P_i$
are said to belong to subdomain $P_i$. Furthermore, triangle pairs, $T_n^{\pm}$,
to an edge in $P_i$, also belong to subdomain $P_i$. Thus, subdomains can share
the same triangle, as seen in Fig.~\ref{fig:meshPartition1}.

Prior to any partitioning, all edges of the mesh are given a
unique edge index, ranging from $1$ to $N$, see \eqref{eq:currentSum}. With a
mesh partitioned into $P$ subdomains, there are $P$ subsets of edge indices,
\begin{equation}
  \mathcal{E}_1,  \mathcal{E}_2, \dots, \mathcal{E}_P \text{,} \quad  
  \bigcup_{i =1}^{P} \mathcal{E}_i = \{1,2,\dots,N \} \text{,} \quad 
  \bigcap_{i =1}^{P} \mathcal{E}_i = \text{\O} \text{,}
  \label{eq:edgeSubset}
\end{equation}
where $\mathcal{E}_i$ are the edge indices of subdomain $i$.

Similarly, all triangles of the mesh are given a unique triangle index, ranging
from $1$ to $N^T$. Here, we use $\mathcal{T}$ to denote all the triangle indices
of the mesh, and $\mathcal{T}_i$ to denote all the triangle indices of subdomain
$i$. 

Consider a triangular mesh partitioned into $P$ subdomains. In order to
calculate elements in the submatrix $Z_{ij}$, $i,j \in [1,P]$,
\eqref{eq:MoMimpedanceMatrix} is calculated with $m \in \mathcal{E}_i$ and $n
\in \mathcal{E}_j$. Upon calculating all the necessary inter-subdomain
calculations (the upper block triangular), the full \gls{MoM} impedance matrix
can be constructed (here denoted $Z^{\mathrm{part}}$). Note that the rows and
columns of the obtained impedance matrix, $Z^{\mathrm{part}}$, are indexed
according to the order of index in $\{\mathcal{E}_1,\dots,\mathcal{E}_P\}$, and
not necessarily to the standard indexing $\{1,\dots,N\}$, \ie, 
\begin{equation}
  Z^{\mathrm{part}} = 
  \kbordermatrix{\mbox{indices}& \mathcal{E}_1 & \mathcal{E}_2 & \dots &
    \mathcal{E}_P\\
    \mathcal{E}_1 & Z_{11} & Z_{12} & \dots & Z_{1P}\\
    \mathcal{E}_2 & Z_{21} & Z_{22} & \dots & Z_{2P}\\
    \vdots & \vdots & \vdots & \ddots & \vdots\\
    \mathcal{E}_P & Z_{P1} & Z_{P2} & \dots & Z_{PP}\\
  } \text{.}
  \label{eq:scrambledMatrix}
\end{equation}
The translation from $Z^{\mathrm{part}}$ to $Z$ is carried out by creating the
translation index $\mathcal{P} = [ \mathcal{E}_1,\dots,\mathcal{E}_P]$, and
letting $Z_{\mathcal{P}(m),\mathcal{P}(n)} =  Z^{\mathrm{part}}_{mn}$.

The integrals \eqref{eq:MoMimpedanceMatrix} are calculated numerically with,
\eg, Gaussian quadrature. The weakly singular integrals arising when $m$ and $n$
in \eqref{eq:MoMimpedanceMatrix} share any supporting triangle are handled with
the \gls{DEMCEM} package \cite{polimeridis_direct_2008,
polimeridis_complete_2010, polimeridis_fast_2011, polimeridis_direct_2011}. By
intersecting the triangle indices between subdomains, the shared triangles are
obtained. When to use Gaussian quadrature or the \gls{DEMCEM} package thus
depends on whether the edges $m$ and $n$ share any supporting triangle.
\gls{DEMCEM} is capable of solving weakly singular integrals of adjacent
triangles, here, however, we only use it on the self-term.

In order to facilitate the numerical calculation of
\eqref{eq:MoMimpedanceMatrix} when $m \in \mathcal{E}_i$, and $n \in
\mathcal{E}_j$ (\ie, the calculation of the submatrix $Z_{ij}$), each subdomain
has the following \emph{designated data} on which edges and triangles belong to
that subdomain:
\begin{equation}
    \begin{cases}
      \mathcal{E}_i \text{,} \quad \text{Global edge indices of subdomain $P_i$} \\
      N_i \text{,} \quad \text{Number of edges in subdomain $P_i$}\\
      \mathcal{T}_i^+ = \mathcal{T}^+(\mathcal{E}_i) \text{,} \text{ } \text{ }
      \text{Global
      plus triangle indices of $P_i$} \\
      \mathcal{T}_i^-  = \mathcal{T}^-(\mathcal{E}_i) \text{,} \text{ }\text{}
      \text{Global minus triangle indices of $P_i$} \\
      \mathcal{T}_i = \mathcal{T}^+_i \cup
      \mathcal{T}^-_i  \text{,} \quad \text{Global triangle
      indices of $P_i$} \\
      N^{T}_i \text{,} \quad \text{Number of triangles in subdomain $P_i$}\\
      \widehat{\mathcal{T}}_i^{+} \text{,} \quad \text{Intrinsic plus triangle
      indices of subdomain $P_i$} \\
      \widehat{\mathcal{T}}_i^{-} \text{,} \quad \text{Intrinsic minus triangle
      indices of subdomain $P_i$} \\
      \widehat{\mathcal{E}}_i = \{1,\dots,N_i\} \text{,} \quad \text{Intrinsic edge
      indices of $P_i$} \\
      \widehat{\mathcal{T}}_i = \{1,\dots,N^T_i\} \text{,} \quad \text{Intrinsic
      triangle indices of $P_i$}
    \end{cases}
  \label{eq:partData}
\end{equation}
where $\mathcal{T}^+$ and $\mathcal{T}^-$ are the global plus and minus triangle
indices, respectively (the $m$th index of $\mathcal{T}^+$ yields the plus
triangle index to the plus triangle $T^+_m$, corresponding to edge $m$). The
intrinsic plus and minus triangle indices relate the intrinsic edge indices
with the intrinsic triangle indices, \eg, $\widehat{\mathcal{T}}^+_i(1)$ yields
the plus triangle index, in the range $[1,N_i^T]$, to the first edge in subdomain $P_i$. They fulfill
\begin{equation}
  \mathcal{T}_i(\widehat{\mathcal{T}}^{\pm}_i(k) ) \implies T_{m}^{\pm} \text{,} \quad
  \mathcal{E}_i(k) = m \text{,} \quad k = 1,\dots,N_i \text{.}
  \label{eq:intrinsicPlusMinus}
\end{equation}
To determine the intrinsic plus and minus triangle indices we
propose Algorithm \ref{algo:plusminusTriangle}. 
 \begin{algorithm}
 \caption{Calculation of the intrinsic plus or minus triangle indices of
 subdomain $i$}
 \begin{algorithmic}[1]
 \renewcommand{\algorithmicrequire}{\textbf{Input:}}
 \renewcommand{\algorithmicensure}{\textbf{Output:}}
 \REQUIRE $\mathcal{T}_i$, $\mathcal{T}_i^{\pm}$ \\
 \ENSURE  $\widehat{\mathcal{T}}_i^{\pm}$\\
 \STATE  $\widehat{\mathcal{T}}_i^{\pm} = $ \texttt{zeros}(1,$N_i$)
 \STATE $\ell = 1$
  \FOR {$j = 1,\dots,N^T_i$}
  \STATE tempIndx $\leftarrow$  find($\mathcal{T}_i^{\pm} = \mathcal{T}_i(j)$)
  \STATE $\widehat{\mathcal{T}}_i^{\pm}(\text{tempIndx}) = \ell$
  \STATE $\ell \leftarrow \ell +1$
  \ENDFOR
 \RETURN $\widehat{\mathcal{T}}_i^{\pm}$
\end{algorithmic} 
\label{algo:plusminusTriangle}
\end{algorithm}

With the subdomain designated data of \eqref{eq:partData} according to two
subdomains, subdomain $P_i$ and $P_j$, the submatrix $Z_{ij} \in \mathbb{C}^{N_i
\times N_j}$ may be calculated using \eqref{eq:MoMimpedanceMatrix} with $m\in
[1,N_i]$, $n\in[1,N_j]$, and by using the intrinsic plus and minus triangle
indices, $\widehat{\mathcal{T}}_i^{\pm}$, and $\widehat{\mathcal{T}}_j^{\pm}$.
Calculating the \gls{MoM} impedance matrix $Z^{\mathrm{part}}$ using the
proposed partitioning strategy is detailed in Algorithm
\ref{algo:partitioningAlgo}.

 \begin{algorithm}
 \caption{Calculation of the \gls{MoM} impedance matrix
 $Z^{\mathrm{part}}$}
 \begin{algorithmic}[1]
 \renewcommand{\algorithmicrequire}{\textbf{Input:}}
 \renewcommand{\algorithmicensure}{\textbf{Output:}}
 \REQUIRE Triangular mesh, $\mathcal{E}$, $\mathcal{T}$,
 \ENSURE $Z^{\mathrm{part}}$, translation index $\mathcal{P}$
 \STATE Partition the mesh ($\mathcal{E}$, $\mathcal{T}$) into $P$ subdomains:
 $S_1,\dots,S_P$\\
 \STATE Allocate $Z^{\mathrm{part}}$\\
 \STATE Calculate the designated data \eqref{eq:partData} for all
 $\{S_n\}_{n=1}^{P}$
  \FOR {$j = 1,\dots,P$}
  \STATE Load designated data \eqref{eq:partData} of subdomain $S_j$
  \FOR {$i = 1,\dots,j$}
  \STATE Load designated data \eqref{eq:partData} of subdomain $S_i$
  \IF{Exist shared triangles: $\mathcal{T}_j \cap \mathcal{T}_i \neq  \text{\O}$}
  \STATE Calculate intrinsic edge indices corresponding to shared triangles in
  $S_i, S_j$:
  $\rightarrow \mathcal{I},\mathcal{J}$
  \STATE Solve \eqref{eq:MoMimpedanceMatrix} using DEMCEM, $m\in \mathcal{I}$,
  $n\in \mathcal{J}$, $\widehat{\mathcal{T}}_i^{\pm}$, and
  $\widehat{\mathcal{T}}_j^{\pm}$ to obtain $Z_{ij}'$
\STATE Solve \eqref{eq:MoMimpedanceMatrix} numerically with,
$m\in [1,N_i]\backslash\mathcal{I}$, $n\in [1,N_j]\backslash\mathcal{J}$,
$\widehat{\mathcal{T}}_i^{\pm}$, and $\widehat{\mathcal{T}}_j^{\pm}$ to obtain
$Z_{ij}^{\prime\prime}$
\STATE $Z_{ij} \leftarrow Z_{ij}' + Z_{ij}^{\prime\prime}$
\ELSE
  \STATE Solve \eqref{eq:MoMimpedanceMatrix} numerically with,
  $m\in [1,N_i]$,
$n\in[1,N_j]$,
$\widehat{\mathcal{T}}_i^{\pm}$, and $\widehat{\mathcal{T}}_j^{\pm}$ to obtain
$Z_{ij}$
  \ENDIF
  \STATE $Z^{\mathrm{part}}_{ij} \leftarrow Z_{ij}$
  \ENDFOR
  \ENDFOR
  \STATE Fill lower triangular part: $Z^{\mathrm{part}}_{ji}
  \leftarrow (Z^{\mathrm{part}}_{ij})^T$, $i \neq j$
  \STATE Determine translation index: $\mathcal{P} =
  \{\mathcal{E}_1,\dots,\mathcal{E}_P\}$
  \RETURN $Z^{\mathrm{part}}$, $\mathcal{P}$
\end{algorithmic} 
\label{algo:partitioningAlgo}
\end{algorithm}

\subsection{Array Partitioning}
In this section, we describe the proposed partitioning algorithm for the array case.
For simplicity, consider elements in a Cartesian $xy$-grid with finite
translations in two dimensions. The elements can be electrically connected and
placed above a finite ground plane.

\begin{figure}[ht!]
  \begin{center}
    \includegraphics[width=\columnwidth]{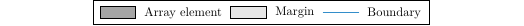}

    \subfloat[]{\includegraphics[width=120pt]{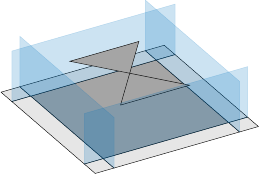}}
    \hfill
    \subfloat[]{\includegraphics[width=120pt]{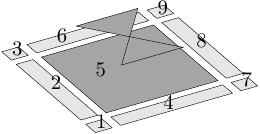}}
  \end{center}
  \caption{Array element and margin. (a) Partitioning boundary highlighted. (b)
  Exploded view with the nine components numbered.}
  \label{fig:UnitCellMargin}
\end{figure}

\begin{figure}[ht!]
  \begin{center}
    \includegraphics[width=\columnwidth]{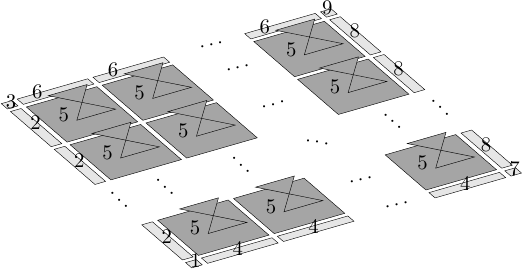}
  \end{center}
  \caption{Complete array antenna constructed from the nine components.}\label{fig:completeARray}
\end{figure}

We require that an array can be built by an array element and a \emph{margin},
as depicted in Fig.~\ref{fig:UnitCellMargin}(a). The array element and margin
are partitioned according to the highlighted boundaries of
Fig.~\ref{fig:UnitCellMargin}(a), which yields the nine subdomains displayed in
Fig.~\ref{fig:UnitCellMargin}(b), from which an $N_x \times N_y$ array may be
constructed, as proposed in Fig.~\ref{fig:completeARray}.

It is not necessary to mesh the complete array of Fig.~\ref{fig:completeARray},
instead; the structure in Fig.~\ref{fig:UnitCellMargin}(a) is meshed and
partitioned into the nine subdomains of Fig.~\ref{fig:UnitCellMargin}(b). The
complete array can then be described by finite translations of the nine
subdomains to build the array in Fig.~\ref{fig:completeARray}. Here, we require
the mesh across the internal borders to be matched by assuming grid points along
the edges periodically in both translation directions.

The nine subdomains of Fig.~\ref{fig:UnitCellMargin}(a) are here denoted
\emph{array components}, whereas a \emph{part} refers to an array component
together with its offset in space. Representing a complete $N_x\times N_y$
electrically connected array (\eg, Fig.~\ref{fig:completeARray}) thus
requires (up to) nine array components with designated data \eqref{eq:partData}
and an offset matrix:
\begin{equation}
  R_{\mathrm{offset}} = 
  \begin{bmatrix}
    x_1 & \cdots & x_{P_N} \\
    y_1 & \cdots & y_{P_N} \\
    z_1 & \cdots & z_{P_N} \\
  \end{bmatrix}
  \text{,} \quad P_N = N_xN_y + 2N_x + 2N_y + 4 \text{,}
  \label{eq:offsetMatrix}
\end{equation}
where $P_N$ is the total number of parts in the array, and $x_i,y_i,z_i$ refer
to the offset coordinates of part $P_i$. The ordering of the part offsets in
$R_{\mathrm{offset}}$ is crucial for obtaining a multilevel block Toeplitz
\gls{MoM} impedance matrix and is further explained in Section
\ref{sec:Toeplitz}.

\begin{figure}[ht!]
  \begin{center}
    \includegraphics[width=\columnwidth]{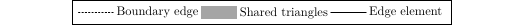}

    \subfloat[]{\includegraphics[width=125pt]{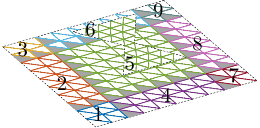}}
    \hfill
    \subfloat[]{\includegraphics[width=125pt]{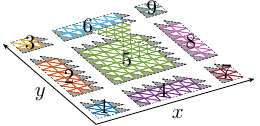}}
  \end{center}
  \caption{Meshed array element and margin configuration. (a) The numbered array
  components (b) Exploded view.}\label{fig:ucmMeshed}
\end{figure}

Since two antenna elements have the same designated data~\eqref{eq:partData},
Algorithm \ref{algo:partitioningAlgo} has to be extended in order to correctly
calculate the shared triangles between two parts of the same array component.
Consider the example of the nine array components meshed in
Fig.~\ref{fig:ucmMeshed}(a), where the shared triangles are highlighted in the
exploded view of Fig.~\ref{fig:ucmMeshed}(b). To determine the shared triangles
between two parts, $P_i$ and $P_j$, $i\neq j$, Algorithm \ref{algo:partitioningAlgo}
performs the intersection calculation $\mathcal{T}_i \cap \mathcal{T}_j$. 
For two adjacent parts of the \emph{same} array component, the intersection
calculation yields all the triangles of the array component. Subsequently,
Algorithm \ref{algo:partitioningAlgo} must be modified such that the shared
triangles calculation of two adjacent parts of the same array component becomes
correct (for two parts of \emph{different} array components, the intersection
calculation suffices).

\begin{figure}[ht!]
  \begin{center}
    \includegraphics[width=\columnwidth]{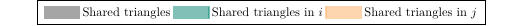}

    \subfloat[]{\includegraphics[width=125pt]{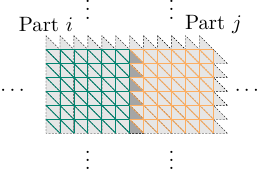}}
    \hfill
    \subfloat[]{\includegraphics[width=125pt]{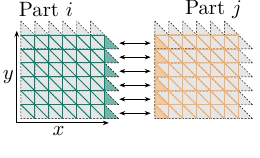}}
  \end{center}
  \caption{Shared triangles between two adjacent identical array components. (a)
  Shared triangles. (b) Shared triangles in each respective
part.}\label{fig:sharedTriangles}
\end{figure}

To address this issue, consider the example in
Fig.~\ref{fig:sharedTriangles}(a), where two adjacent parts of the same array
component, parts $i$ and $j$, are displayed. The shared triangles between the two
parts are located along their connection and are highlighted in
Fig.~\ref{fig:sharedTriangles}(b). Here, it may be noted that the rightmost
shared triangles of part $i$ correspond to the leftmost shared triangles of
part $j$, \ie, in Fig.~\ref{fig:sharedTriangles}(b) the topmost blue triangle
corresponds to the topmost red triangle, \etc. Thus, to determine which shared
triangle in part $i$ corresponds to which shared triangle in part $j$, the
shared triangles on the right side of part $i$ and on the left side of part $j$
must be determined. Furthermore, these two sets of shared triangles must be
sorted according to their $y$-position, to match their correspondence in the other
part (for two parts located over each other, the shared triangles of the top and
bottom sides are sorted according to their $x$-position. In the 3D case, the
$z$-position may be used to sort triangles with the same $x$\slash$y$-position).

By identifying all shared triangles of each array component [see
Fig.~\ref{fig:ucmMeshed}(b)], determining on which side they are located (left,
right, top, bottom, south corner, and north corner), and then sorting them
according to their position, corresponding shared triangles in two parts of the
same array component can be calculated. Calculating the shared triangles on a
side of an array component can be done by intersecting its shared triangles with
the shared triangles of another array component that is located to the
left$\slash$right or on top$\slash$bottom of itself, \eg, the left shared
triangles of array component \num{5} are obtained intersecting the shared
triangles of array component \num{5} and \num{2}.
The proposed procedure leads to a matrix of indices to store the shared
triangle indices of each array component, for each side:
\begin{equation}
  \mathcal{T}^{\mathrm{Shared}} = 
\begin{bmatrix}
  & \mathcal{T}^{r1} & \mathcal{T}^{t1} & & &\\
  & \mathcal{T}^{r2} & \mathcal{T}^{t2} & \mathcal{T}^{b2}  & \mathcal{T}^{sc2}& \\
  & \mathcal{T}^{r3} & & \mathcal{T}^{b3}  & \mathcal{T}^{sc3}& \\
  \mathcal{T}^{l4} & \mathcal{T}^{r4} & \mathcal{T}^{t4} & & & \mathcal{T}^{nc4}\\
  \mathcal{T}^{l5} & \mathcal{T}^{r5} & \mathcal{T}^{t5} &\mathcal{T}^{b5} &
  \mathcal{T}^{sc5}&\mathcal{T}^{nc5} \\
  \mathcal{T}^{l6} & \mathcal{T}^{r6} &  &\mathcal{T}^{b6} & \mathcal{T}^{sc6}& \\
  \mathcal{T}^{l7} &  & \mathcal{T}^{t7} & & & \mathcal{T}^{nc7}\\
  \mathcal{T}^{l8} &  & \mathcal{T}^{t8} & \mathcal{T}^{b8}& &  \mathcal{T}^{nc8}\\
  \mathcal{T}^{l9} &  & & \mathcal{T}^{b9} & & \\
\end{bmatrix}
\text{,}
  \label{eq:sharedTrianglesParttype}
\end{equation}
where $\mathcal{T}^{si}$ is the indices of shared triangles of array component
$i$, on the left ($s=l$), right ($s=r$), top ($s=t$), bottom side ($s=b$), south
corner ($s = sc$), or north corner ($s = nc$).

\begin{figure}[ht!]
  \begin{center}
    \subfloat[]{\includegraphics[width=84pt]{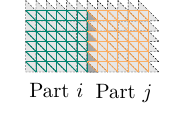}}
    \hfill
    \subfloat[]{\includegraphics[width=84pt]{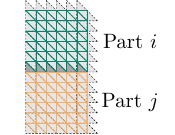}}
    \hfill
    \subfloat[]{\includegraphics[width=84pt]{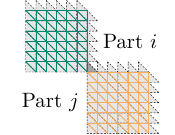}}
  \end{center}
  \caption{Three relative positions of two adjacent parts such that they share
  triangles. (a) Two parts besides each other. (b) Two parts on top of each
other. (c) Two parts placed on the north west diagonal.}\label{fig:sharedTrianglesOptions}
\end{figure}

There are three scenarios in which two adjacent parts of any array component
may share triangles. The three scenarios are displayed in
Fig.~\ref{fig:sharedTrianglesOptions}.

Let the notation $\mathcal{T}^{\mathrm{Shared}}_{ij}$ indicate the $ij$th
submatrix of $\mathcal{T}^{\mathrm{Shared}}$ \eqref{eq:sharedTrianglesParttype},
and consider the two parts $i$ and $j$ to be of array component $A$ and $B$,
respectively, where $A,B\in[1,9]$. The shared triangles in the two parts, for
each of the three scenarios in Fig.~\ref{fig:sharedTrianglesOptions} are then as
follows:
%
%
\begin{enumerate}[label=(\alph*)]
  \item Shared triangles in part $i$: $\mathcal{T}^{\mathrm{Shared}}_{A,2}$, shared
    triangles in part $j$: $\mathcal{T}^{\mathrm{Shared}}_{B,1}$.

  \item Shared triangles in part $i$: $\mathcal{T}^{\mathrm{Shared}}_{A,4}$, shared
    triangles in part $j$: $\mathcal{T}^{\mathrm{Shared}}_{B,3}$.

  \item Shared triangles in part $i$: $\mathcal{T}^{\mathrm{Shared}}_{A,5}$,
    shared triangles in part $j$: $\mathcal{T}^{\mathrm{Shared}}_{B,6}$.
\end{enumerate}
Note that depending on the mesh, multiple corners may or may not overlap.

Each part of the array (consisting of an array component and an offset in
space), has a row and column location, as indicated in
Fig.~\ref{fig:completeARray}. To calculate the submatrix $Z_{ij}$, according to
parts $i$ and $j$, the parts' array components, designated data \eqref{eq:partData},
and row$\slash$column location are loaded. The parts' respective row and column
location govern if they are situated in any of the three configurations
displayed in Fig.~\ref{fig:sharedTrianglesOptions}, thus deciding if any
shared triangle calculation is necessary [and if so, which case of (a), (b), and
(c)].

The row$\slash$column location of each part in Fig.~\ref{fig:completeARray}
yields a $3\times P_N$ matrix (similar to the offset matrix), with the row,
column, and array component data of each part:
\begin{equation}
  M = \begin{bmatrix}
    \mathrm{row}_1 & \dots & \mathrm{row}_{P_N}\\
    \mathrm{col}_1 & \dots & \mathrm{col}_{P_N}\\
    \mathrm{comp}_1 & \dots & \mathrm{comp}_{P_N}\\
  \end{bmatrix} \in \mathbb{Z}^{3\times P_N} \text{.}
  \label{eq:RowColumn}
\end{equation}
where $\mathrm{comp}_i \in [1,9]$ denotes the array component of part $i$,
and where each column index of matrix $M$ corresponds to a part. Subsequently, each
part has an 'identity', ranging from 1 to $P_N$, where $P_N$ is given in
\eqref{eq:offsetMatrix}. 

With the nine array components, offset and row$\slash$column matrices, and a
method for handling shared triangles between parts of the same array component,
the \gls{MoM} impedance matrix $Z^{\mathrm{part}}$ can be obtained. Calculating
all submatrices $Z_{ij}$, $i,j \in [1,P_N]$, is, however, not necessary. In
Section \ref{sec:Toeplitz}, the finite translation symmetry of the array is
exploited to yield a reduction in the number of submatrix calculations to
obtain the complete \gls{MoM} matrix $Z^{\mathrm{part}}$.

\subsection{Multilevel Toeplitz Impedance Matrix Representation}
\label{sec:Toeplitz}
\begin{figure}[ht!]
  \begin{center}
    \includegraphics[width=\columnwidth]{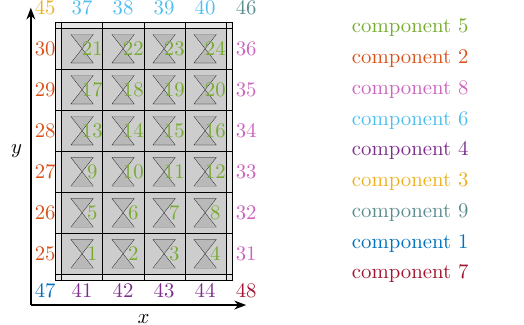}
  \end{center}
  \caption{Indexing order of the parts, where the array component of each part
  is highlighted by the color of the indexing number.}\label{fig:partsIndexing}
\end{figure}

Since the Greens function $G(\bm{r},\bm{r}')$ in \eqref{eq:MoMimpedanceMatrix}
is translation invariant, a submatrix calculation between two parts, $i$ and
$j$, yields the same submatrix as for the case where both parts are displaced by
the same amount in space. For example, consider the $4\times6$ array depicted in
Fig.~\ref{fig:partsIndexing}. The submatrix calculation between parts 1 and 6
yields the same submatrix obtained for parts 2 and 7. Similarly, parts
(25,5) and (27,13) yield the same submatrix.

Upon identifying all the unique edge element distance relations between all
parts, only the necessary submatrix calculations are obtained. Additionally,
calculating the submatrices in a specific order (\eg, calculating submatrices
corresponding to parts of array component 5 in Fig.~\ref{fig:ucmMeshed} first), the resulting \gls{MoM}
impedance matrix $Z^{\mathrm{part}}$ inherits a multilevel block Toeplitz
structure\cite{brandt-moller_extended_2024, helander_synthesis_2018,
kindt_array_2003}. 

Recall the matrix $M$ containing the row$\slash$column and array component data
according to each part. Each part is given an index, ranging from $1$ to $P_N$.
The index of a part depends on its array component.  Here, parts of array
component \num{5} (the antenna element) are indexed first, where the complete
indexing order of the array components is given as $5, 2,8,6,4,3,9,1,7$.
Furthermore, of all parts of the same array component, the part leftmost and
closest to the bottom is indexed first, with the second being the part to the
right of the first part. In Fig.~\ref{fig:partsIndexing} each part of the
$4\times 6$ array is indexed according to the proposed indexing strategy.

The resulting \gls{MoM} impedance matrix is obtained on the form
\begin{equation}
  Z^{\mathrm{part}} = 
  \begin{bmatrix}
    A & B^T \\
    B & C
  \end{bmatrix}
  \text{,}
  \label{eq:matrixForm}
\end{equation}
where matrix $A$ corresponds to the interaction between the array element parts,
matrix $B$ corresponds to the margin-to-array elements interaction, and matrix $C$
corresponds to the interaction between the margin parts.

As already established, only a fraction of submatrix calculations are necessary to
construct matrices $A$, $B$, and $C$. Starting with matrix $A$, the largest
of the three matrices, its structure is obtained as a multilevel block Toeplitz
matrix on the form 
\begin{equation}
  \begin{aligned}
    A = A^{[2]} = 
  & \begin{bmatrix}
    A_{0}^{[1]}  & \cdots &  A_{N_y-1}^{[1]T}\\
    \vdots & \ddots & \vdots \\
    A_{N_y-1}^{[1]}  & \cdots &  A_{0}^{[1]}\\
  \end{bmatrix}
  \text{,}
  \\
  & A_{i}^{[1]}  = 
  \begin{bmatrix}
    A_{0}^{[0]}(i)  & \cdots &  A_{1-N_x}^{[0]}(i)\\
    \vdots & \ddots & \vdots \\
    A_{N_x-1}^{[0]}(i)  & \cdots &  A_{0}^{[0]}(i)\\
  \end{bmatrix}
  \text{,} 
  \\
  &\quad \quad A_{j}^{[0]} (i) = Z_{iN_x + 1 + (|j|+j)/2,1+ (|j| - j)/2}
  \text{.}
\end{aligned}
  \label{eq:submatrixAA}
\end{equation}
Here, the notation $A^{[L]}$ denotes an $L$-level block Toeplitz matrix
\cite{choromanski2022block,MoMHH}.

\begin{figure}[ht!]
  \begin{center}
    \subfloat[]{\includegraphics[width=115pt]{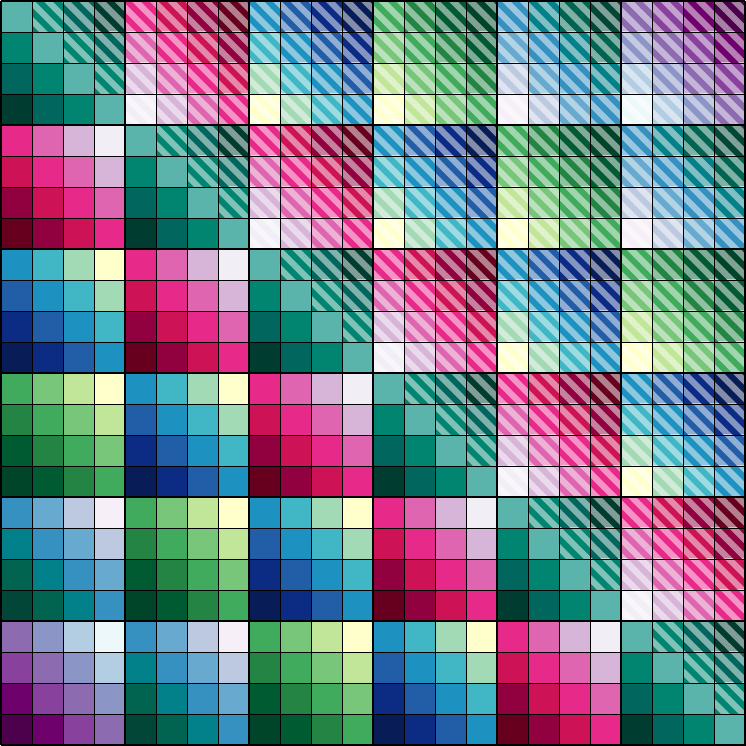}}
    \hfill
    \subfloat[]{\includegraphics[width=115pt]{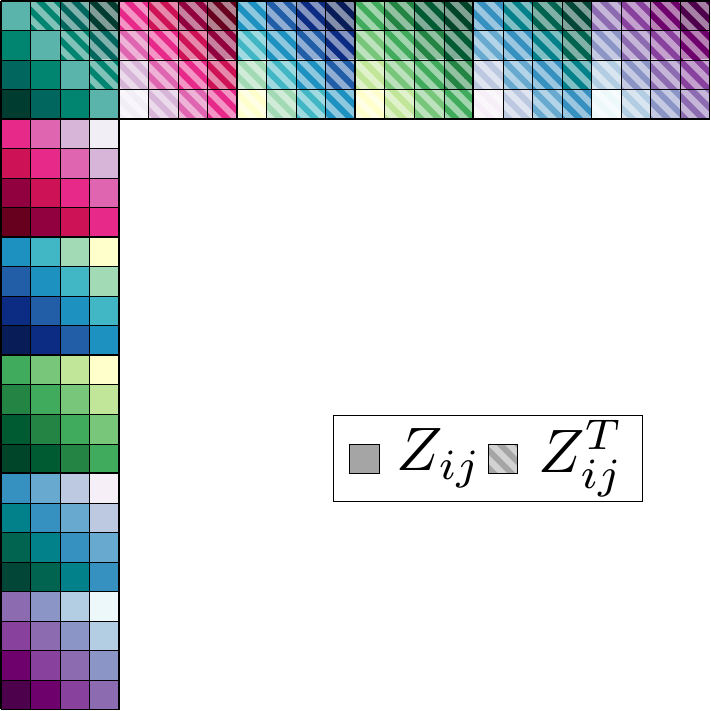}}
    \hfill
    \subfloat[]{\includegraphics[width=19.4pt]{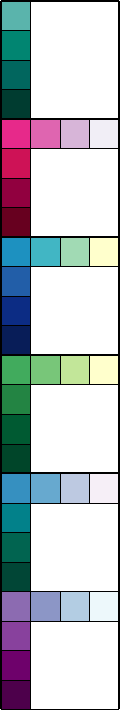}}
  \end{center}
  \caption{Visualization of matrix $A$ for the $4\times 6$ array case. Blocks of
  the same color indicate identical submatrices. (a) The complete matrix $A$.
  (b) The first level of Toeplitz structure of $A$. (c) The second level of
  Toeplitz structure, and the symmetry of $A$ exploited.}\label{fig:Amatrices}
\end{figure}

In Fig.~\ref{fig:Amatrices}(a), a visualization of matrix $A$ is depicted. The
displayed matrix structure corresponds to the $4 \times 6$ array depicted in
Fig.~\ref{fig:partsIndexing}. Exploiting the first level of block Toeplitz
structure reduces the needed calculations to Fig.~\ref{fig:Amatrices}(b). By
further exploiting the second level of block Toeplitz structure (and
symmetry), the unique submatrix calculations in Fig.~\ref{fig:Amatrices}(c) are
obtained.

Thus, for an $N_x \times N_y$ array, the $A$ matrix requires $N_x N_y +
(N_x-1)(N_y-1)$ submatrix calculations.

Matrix $B$, corresponding to the impedance relation between the margin parts and
the array elements, is obtained on the form
\begin{equation}
  \begin{aligned}
  B = 
  \begin{bmatrix}
    B_1^T & B_2^T &  B_3^T & B_4^T &  B_5^T & B_6^T &  B_7^T & B_8^T 
  \end{bmatrix}^T
  \text{.}
\end{aligned}
  \label{eq:submatrixBB}
\end{equation}
Toeplitz structure is only found in matrices $B_1$, $B_2$, $B_3$, and
$B_4$. Matrices $B_5$, $B_6$, $B_7$, and,
$B_8$, corresponding to the impedance relation between the four margin corners
and the array elements, require full submatrix calculations.

Matrices $B_1$ and $B_2$ inherit the same block Toeplitz \emph{structure}. This
block Toeplitz structure is obtained as 
\begin{equation}
  \begin{aligned}
    B_i = B_i^{[1]}& = 
  \begin{bmatrix}
    B^{[0]}_{0}(i) & \cdots & B^{[0]}_{1-N_y}(i)\\
    \vdots & \ddots & \vdots \\
    B^{[0]}_{N_y-1}(i) & \cdots & B^{[0]}_{0}(i)\\
  \end{bmatrix} \text{,}\\
    B_{j}^{[0]}(i) &= 
  \begin{bmatrix}
    Z_{p,q} & \cdots & Z_{p,q+N_x}
  \end{bmatrix}\text{,} \\
    p & = N_xN_y +(i-1)N_y+ (|j|+j)/2+1 \text{,} \\
      &q = [(|j|-j)/2] N_x+1 \text{,}
\end{aligned}
  \label{eq:submatrixB12}
\end{equation}
where $i = 1,2$.

Similarly, matrices $B_3$ and $B_4$ inherit the same block Toeplitz
structure as
\begin{equation}
  \begin{aligned}
    B_i =&
  \begin{bmatrix}
    B^{[1]}_{0}(i) & \cdots & B^{[1]}_{N_y-1}(i)\\
  \end{bmatrix} \text{,}\\
    B_{j}^{[1]}&(i) =
  \begin{bmatrix}
    B^{[0]}_{1}(i,j) & \cdots & B^{[0]}_{1-N_x}(i,j) \\
    \vdots & \ddots & \vdots \\
    B^{[0]}_{N_x-1}(i,j) & \cdots & B^{[0]}_{1}(i,j) \\
  \end{bmatrix} \text{,} \\
    B^{[0]}_{k}&(i,j) = Z_{p,q} \text{,}\\
    p & = N_xN_y +2N_y + (i-3)N_x + (|k|+k)/2+1 \text{,} \\
      &q = jN_x +1 + (|k|-k)/2 \text{,}
\end{aligned}
  \label{eq:submatrixB34}
\end{equation}
where $i = 3,4$.

\begin{figure}[ht!]
  \begin{center}
    \subfloat[]{\includegraphics[width=132.35pt]{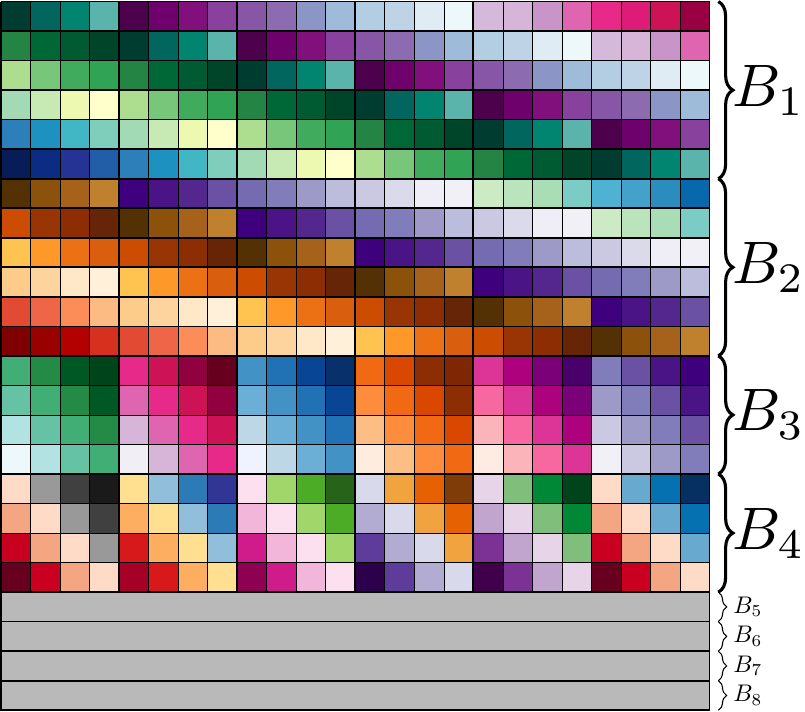}}
    \hfill
    \subfloat[]{\includegraphics[width=117.65pt]{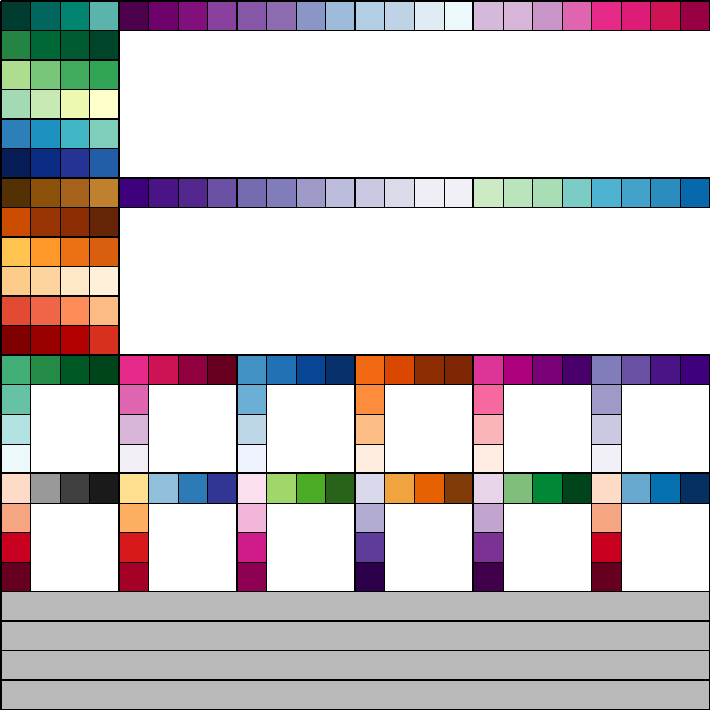}}
  \end{center}
  \caption{Visualization of matrix $B$ for the $4\times 6$ array case. Blocks
  of the same color indicate identical submatrices. Grey areas indicate unique
  submatrices. (a) Complete matrix $B$. (b) First level of Toeplitz
  structure.}\label{fig:Bmatrices}
\end{figure}

In Fig.~\ref{fig:Bmatrices}(a), a visualization of matrix $B$ is depicted. The
displayed matrix structure corresponds to the $4\times6$ array shown in
Fig.~\ref{fig:partsIndexing}. Exploiting the only level of block Toeplitz
structure available in matrix $B$ yields the reduction of the number of
submatrices of Fig.~\ref{fig:Bmatrices}(b). The grey parts in
Fig.~\ref{fig:Bmatrices}(a) and Fig.~\ref{fig:Bmatrices}(b) correspond to
matrices $B_i$, $i=5,\dots,8$, and require full submatrix
calculations.

Matrix $C$, corresponding to the impedance relation among the margin parts, is
obtained on the form
\begin{equation}
  \begin{aligned}
  C = 
  \begin{bmatrix}
    C_{11} & \cdots & C_{15}\\
    \vdots & \ddots & \vdots\\
    C_{51} & \cdots & C_{55}\\
  \end{bmatrix}
  \text{,}
\end{aligned}
  \label{eq:submatrixCC}
\end{equation}
where $C_{ij} = C_{ji}^T$. Block Toeplitz structure is only found in matrices
$C_{11}$, $C_{21}$, $C_{22}$, $C_{33}$, $C_{43}$, and $C_{44}$. 

Matrices $C_{11}$, $C_{21}$, and $C_{22}$ inherit a block Toeplitz structure on the
form 
\begin{equation}
  \begin{aligned}
    C_{ij} = C^{[1]}(ij) =&  
  \begin{bmatrix}
    C_0^{[0]}(i,j) & \cdots & C_{1-N_y}^{[0]}(i,j) \\
    \vdots & \ddots & \vdots \\
    C_{N_y-1}^{[0]}(i,j) & \cdots & C_0^{[0]}(i,j) \\
  \end{bmatrix} \text{,} \\
    C_{k}^{[0]}(i,j) &= Z_{p+(|k|+k)/2,q + (|k| - k)/2} \\
    p & = N_xN_y + (i-1)N_y  + 1 \text{,} \\
    q & = N_xN_y + (j-1)N_y + 1\text{,}
\end{aligned}
  \label{eq:submatrixC123}
\end{equation}
where $i,j \in\{ 1,2\}$.

Similarly, matrices $C_{33}$, $C_{43}$, and $C_{44}$ inherit a block Toeplitz
structure on the form 
\begin{equation}
  \begin{aligned}
    C_{ij} = C^{[1]}(ij) =&  
  \begin{bmatrix}
    C_0^{[0]}(i,j) & \cdots & C_{1-N_x}^{[0]}(i,j) \\
    \vdots & \ddots & \vdots \\
    C_{N_x-1}^{[0]}(i,j) & \cdots & C_0^{[0]}(i,j) \\
  \end{bmatrix} \text{,} \\
    C_{k}^{[0]}(i,j) &= Z_{p+(|k|+k)/2,q + (|k| - k)/2} \\
    p & = N_xN_y +2N_y +(i-3)N_x  + 1 \text{,} \\
    q & = N_xN_y + 2N_y+(j-3)N_x + 1\text{,}
\end{aligned}
  \label{eq:submatrixC456}
\end{equation}
where $i,j \in\{ 3,4\}$.

\begin{figure}[ht!]
  \begin{center}
    \subfloat[]{\includegraphics[width=117.7pt]{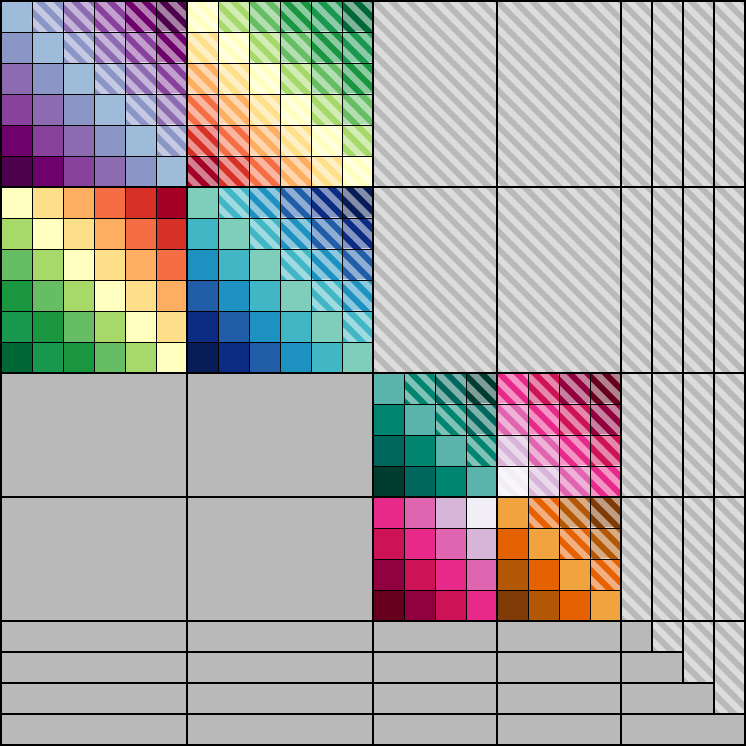}}
    \hfill
    \subfloat[]{\includegraphics[width=117.7pt]{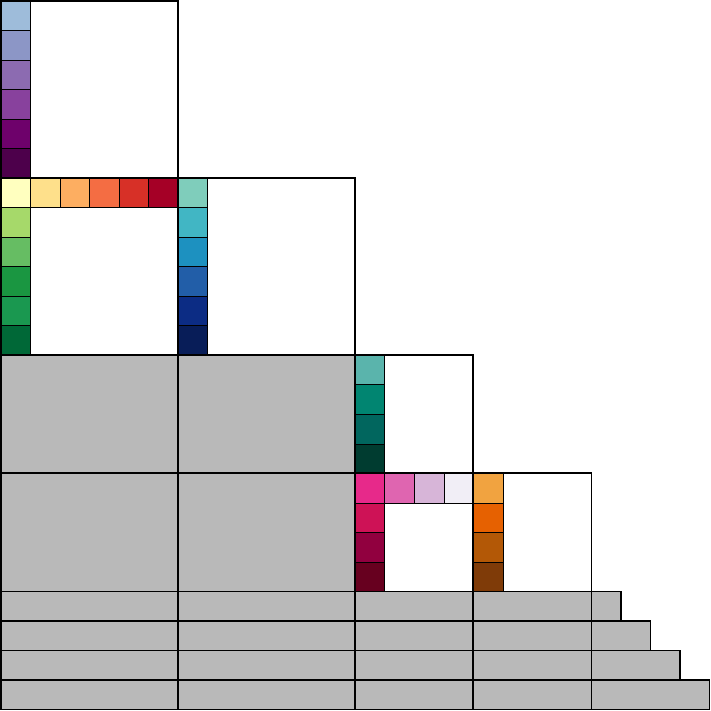}}
  \end{center}
  \caption{Visualization of matrix $C$ for the $4\times 6$ array case. Blocks of
    the same color indicate identical submatrices. Grey areas indicate unique
    submatrices. (a) Complete matrix $C$. (b) First level of Toeplitz
    structure.}\label{fig:Cmatrices}
\end{figure}

In Fig.~\ref{fig:Cmatrices}(a), a visualization of matrix $C$ is depicted. The
displayed matrix structure corresponds to the $4\times6$ array shown in
Fig.~\ref{fig:partsIndexing}. Exploiting the only level of block Toeplitz
available in matrix $C$ yields the necessary calculations of
Fig.~\ref{fig:Cmatrices}(b). The grey parts in Fig.~\ref{fig:Cmatrices}(b)
corresponds to matrices $C_{31}$, $C_{41}$, $C_{51}$, $C_{32}$, $C_{42}$,
$C_{52}$, $C_{53}$, $C_{54}$, and the lower triangular of $C_{55}$, which
require full submatrix calculations.

Exploiting the block Toeplitz structure of matrices $A$, $B$, and $C$, the
memory necessary to represent the \gls{MoM} impedance matrix $Z^{\mathrm{part}}$
is strongly reduced. In the Appendix, the exact memory necessary to allocate for
an $N_x\times N_y$ array with the nine-component representation is
determined. Let $e_i$, $i\in[1,9]$ denote the number of edge
elements in array component $i$.  With the assumption $e_2 \approx e_8 \approx
e_6 \approx e_4 = e_m$, $e_3 \approx e_9 \approx e_1 \approx e_7 = e_c$, and
$e_5 = e_E$, the expression for the necessary memory to allocate, $MEM(N_x,N_y)$,
using \eqref{eq:calcScalingA}, \eqref{eq:calcScalingB}, and
\eqref{eq:calcScalingA} reduces to 
\begin{equation}
  \begin{aligned}
    MEM(N_x,N_y) = \left(2N_x N_y - N_x - N_y + 1\right)e_E^2 &\\
    + (8N_xN_y-2N_x-2N_y)e_Ee_m + 4N_xN_y e_E &e_c\\
    + (4N_xN_y + 4N_x+4_Ny -2)&e_m^2\\
    + 8(N_x +N_y)e_me_c + &10e_c^2 \text{.}
  \end{aligned}
  \label{eq:calcScalingApproxNew}
\end{equation}

Using the proposed method, the memory necessary to allocate the \gls{MoM}
impedance matrix scales as $\mathcal{O}(N_xN_ye_E^2)$, assuming $e_E > e_m, e_c$.

\subsection{Solving the System of Linear Equations}
It is advantageous to solve \eqref{eq:MoMMatrixEq} such that effective memory
storage can be utilized. In our other work \cite{MoMHH}, we present two methods
for solving \eqref{eq:MoMMatrixEq} using the sparse representation of the
\gls{MoM} impedance matrix. The first method is a block matrix-valued Rybicki
algorithm that leverages the Schur complement of $Z$ and exploits the first
level of Toeplitz structure in the matrix $A$, to solve \eqref{eq:MoMMatrixEq},
\ie, 
\begin{equation}
  I_1 = U-FI_2 \text{,} \quad I_2 = -\left(C-BF\right)^{-1}BU \text{,} \quad
  I = \begin{bmatrix}
    I_1 \\
    I_2
  \end{bmatrix} 
  \text{,}
  \label{eq:solveEqSystem1}
\end{equation}
where, $U$ and $F$ are obtained by solving the
equation system
\begin{equation}
  A \begin{bmatrix}
    U & F
  \end{bmatrix} 
  = 
  \begin{bmatrix}
    V_1 &  B^T
  \end{bmatrix} 
  \text{,}
  \quad 
  V =  \begin{bmatrix}
    V_1 \\
    0
  \end{bmatrix} 
  \text{,}
  \label{eq:solveEqSystem2}
\end{equation}
where $V_1 = [\hat{\bm{e}}_{i_1},\hat{\bm{e}}_{i_2},\dots,\hat{\bm{e}}_{i_p}]$,
$\hat{\bm{e}}_i$ is the $i$th column of the identity matrix, $i_1$ is the index
to the feeding edge in the array element component (component 5), and where $i_n
= i_1 + (n-1) e_5$ ($e_5$ is the number of edges in the array component).
Additionally, $p$ is the number of considered ports in the antenna.

The Rybicki algorithm solves \eqref{eq:solveEqSystem2} using the first level of
block Toeplitz structure of matrix $A$. For using only the second level of block
Toeplitz structure of matrix $A$, a \gls{GMRES} can be used \cite{MoMHH}.

\subsection{Calculating the Far-field}
\label{sec:farfield}
In the far-field region $kr' \gg kd$, the scattered electric field $\bm{E}^{s}$
can be expressed in terms of the spatial Fourier transform of the surface
current, \ie,  
\begin{equation}
  \bm{E}^{s}(\bm{r}) = \mathrm{j}k \eta \frac{\mathrm{e}^{\mathrm{j}kr}}{4\pi
    r}\hat{\bm{r}} \times \left( \hat{\bm{r}} \times \int_{S}
  \bm{J}(\bm{r}') \mathrm{e}^{\mathrm{j}k\hat{\bm{r}}\cdot
\bm{r}'}\mathrm{d}S'\right)
  \text{.}
  \label{eq:efieldFF}
\end{equation}
The leading order of the far-field amplitude relates to the scattered electric
field in the far-field region as
\begin{equation}
  \bm{E}^{s}(\bm{r}) = \frac{\mathrm{e}^{-\mathrm{j}kr}}{4\pi r}
  \bm{F}(\hat{\bm{r}}) \text{.}
  \label{eq:farfield}
\end{equation}
Here, the polarization component, $\hat{\bm{e}}_{\tau}^*(\hat{\bm{r}})$,
$\tau = \mathrm{co},\mathrm{cx}$, of $\bm{F}(\hat{\bm{r}})$ is obtained as 
\begin{equation}
  F_{\tau
  }(\hat{\bm{r}}) = 
  -\mathrm{j} k \eta
\hat{\bm{e}}_{\tau}^*(\hat{\bm{r}})
  \cdot \int_{S} \bm{J}(\bm{r}')
  \mathrm{e}^{\mathrm{j}k\hat{\bm{r}}\cdot \bm{r}'}\mathrm{d}S'
  \label{eq:farfieldThetaPhi}
\end{equation}
Substituting the current $\bm{J}(\bm{r})$ in the above equation with
\eqref{eq:currentSum} and rearranging yields 
\begin{equation}
  \begin{aligned}
    F_{\tau
    }(\hat{\bm{r}}) 
    \approx -\mathrm{j}k\eta\sum_{n=1}^{N_i} I_n
    \hat{\bm{e}}_{\tau}^*(\hat{\bm{r}}) \cdot  \int_S\bm{f}_n(\bm{r}')
    \mathrm{e}^{\mathrm{j}k\hat{\bm{r}}\cdot \bm{r}'}\mathrm{d}S' \text{,}
\end{aligned}
  \label{eq:farfieldThetaPhiF}
\end{equation}
where $N_i$ is the number of edges in array component $i$, recall
\eqref{eq:partData}. Exploiting that 
\[
  \mathrm{e}^{\mathrm{j}k\hat{\bm{r}}\cdot (\bm{r}' + \bm{d})} =
  \mathrm{e}^{\mathrm{j}k\hat{\bm{r}}\cdot
  \bm{d}}\mathrm{e}^{\mathrm{j}k\hat{\bm{r}}\cdot \bm{r}'} \text{,}\quad
  \bm{d} \in \mathbb{R}^3\text{,}
\]
the calculation of \eqref{eq:farfieldThetaPhiF} for the complete array is
facilitated, and split into two steps: For the observation directions
$\hat{\bm{r}}_m$, the two polarizations $\tau = \mathrm{co}, \mathrm{cx}$, and
for the $N_i$ basis functions of the array component $i$, the tensor
$\hat{F}$ is defined as 
\begin{equation}
  \begin{aligned}
    \hat{F}_{m,\tau,n}
    = \hat{\bm{e}}_{\tau}^*(\hat{\bm{r}})
    \cdot \int_S\bm{f}_n(\bm{r}')
    \mathrm{e}^{\mathrm{j}k\hat{\bm{r}}_m\cdot
  \bm{r}'}\mathrm{d}S' \text{,}
\end{aligned}
  \label{eq:farfieldBasis}
\end{equation}
The far-field of all parts $P$ of the same array component $i$ is then calculated as
the tensor contraction with the current coefficients (basis coefficients)
$I_{n,p}$ and the spatial offsets $\bm{d}_p$, $p = 1,\dots,P$ of the parts:
\begin{equation}
  F_{m,\tau} = -\mathrm{j}k\eta\sum_{p = 1}^{P} \sum_{n = 1}^{N_i}
  \mathrm{e}^{\mathrm{j}k\hat{\bm{r}}_m\cdot \bm{d}_p} I_{n,p} \hat{F}_{m,\tau,n}
  \label{eq:tensorProd}
\end{equation}
It is thus only necessary to calculate \eqref{eq:farfieldBasis} once for each
array component and then sum up the tensor contractions \eqref{eq:tensorProd}
of each array component to obtain the far-field of the array antenna, using the
nine-component representation.

\subsection{Calculating the Surface Current}
With the proposed array decomposition method, the surface current coefficients
$I$ of large antennas can be calculated efficiently. From these coefficients,
most antenna parameters can be calculated, \eg, the far-field, and the
S-parameters. Additionally, visualizing the surface currents is useful in the
design process of an antenna \cite{behdad_compact_2005}.

Visualizing the total surface current on the array is possible with the proposed
array decomposition method, and it is carried out in two steps. First, the current
is visualized on each triangle of each array part. Then, all array parts are
added together, resulting in a current visualization of the whole array. In the
last step, the currents on the shared triangles among the array parts are added,
resulting in a continuous surface current.

\section{Numerical Results}
\label{sec:Numerical}
In this section, we demonstrate the effectiveness of our proposed array
partitioning algorithm by calculating the scattering parameters and the
far-field of two large arrays. First, a validation example is presented, where
the surface currents are calculated. Then, the magnitude of the
memory scaling of the algorithm is presented.

\subsection{Surface Current of a $2\times 3$ Array}
In this validation example, the surface current of a $2\times 3$ array is
calculated using the proposed array decomposition method. 

\begin{figure}[ht!]
  \begin{center}
    \subfloat[]{\includegraphics[width=138pt]{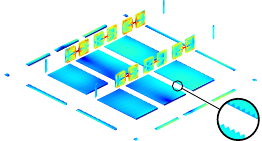}}
    \hfill
    \subfloat[]{\includegraphics[width=114pt]{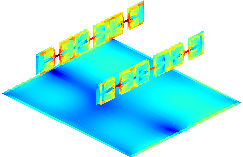}}
  \end{center}
  \caption{Visualized surface current of a $2 \times 3$ array. (a)
  Magnitude of the surface current on each individual array part. (b) Total
  magnitude of the surface current on the area.}
  \label{fig:surfaceCurrent}
\end{figure}

In Fig.~\ref{fig:surfaceCurrent}(a), the magnitude of the surface current of
each individual array part is calculated and visualized. Notice how the
magnitude of the surface current is discontinuous across two parts (see the
highlighted region). This discontinuity is due to the shared triangles not being
added up yet.

In Fig.~\ref{fig:surfaceCurrent}(b), the magnitudes of the surface currents of
the array parts are added together, where the currents on the shared triangles
are summed up. The resulting magnitude is observed to be
continuous across the parts.

\subsection{Scattering Parameters of two Closely Placed $9\times9$ Arrays}
By solving \eqref{eq:MoMMatrixEq} with an excitation voltage matrix $V =
[\hat{\mathbf{e}}_{i_1},\hat{\mathbf{e}}_{i_2},\dots,\hat{\mathbf{e}}_{i_p}]$,
where $p$ is the number of ports of the antenna, and $i_k$, $k=1,\dots,p$ are
the edge indices to the ports, the current matrix $I$ is obtained. Selecting
only the rows according to the indices $i_1,i_2,\dots,i_p$, of the current
matrix $I$ yields the port currents $\hat{I}\in\mathbb{C}^{p\times p}$, from which the admittance matrix
$Y = L\hat{I}L$ is obtained, where $L =
\mathrm{diag}(1/\ell_{i_1},\dots,1/\ell_{i_p})$, and where $\ell_{i_k}$ is the
edge length of edge $i_k$. The inverse of the admittance matrix $Y$ yields the
port impedance matrix $Z_{\mathrm{port}}$, which with a reference impedance
vector $Z_0$, and \cite[Eq.~18]{kurokawa_power_1965} yields the scattering
matrix $S$.

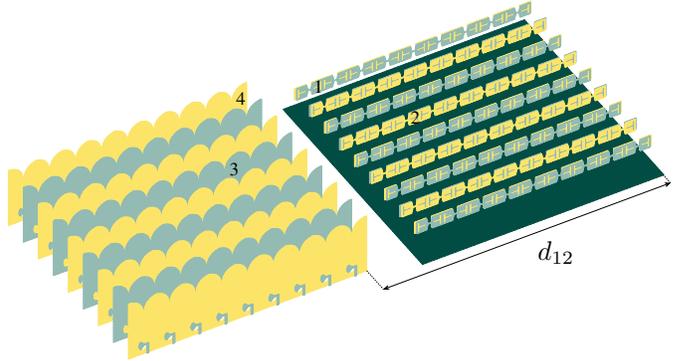
\begin{figure}[ht!]
  \begin{center}
    \begin{tikzpicture}
      \node[inner sep=0pt]  at (0pt,0)
    {\includegraphics
    [width=252pt]{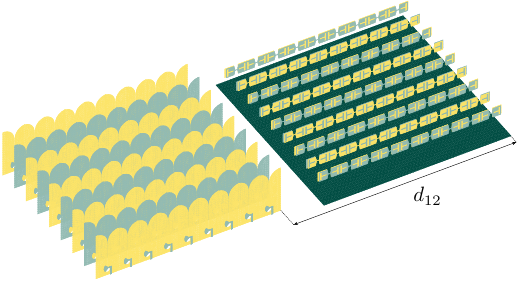}
    };
  \node[scale = 0.7, inner sep=1pt, outer sep=0pt, text = black, line width
    = 0.2pt]
    at (-0.26,1.26) {1};

  \node[scale = 0.7, inner sep=1pt, outer sep=0pt, text = black, line width
    = 0.2pt]
    at (-1.29,1.09) {4};

  \node[scale = 0.7,inner sep=1pt, outer sep=0pt, text = black, line width =
    0.2pt]
    at (1.03,0.84) {2};

  \node[scale = 0.7,inner sep=1pt, outer sep=0pt, text = black, line width =
    0.2pt]
    at (-1.38,0.16) {3};

\end{tikzpicture}
  \end{center}
  \caption{Visualization of two closely placed $9 \times 9$ arrays, where the
  selected ports are highlighted, and where $d_{12} =$ \SI{225.5}{\milli\meter}. }\label{fig:2arraysGeo}
\end{figure}

\begin{figure}[ht!]
  \begin{center}
    \includegraphics[width=\columnwidth]{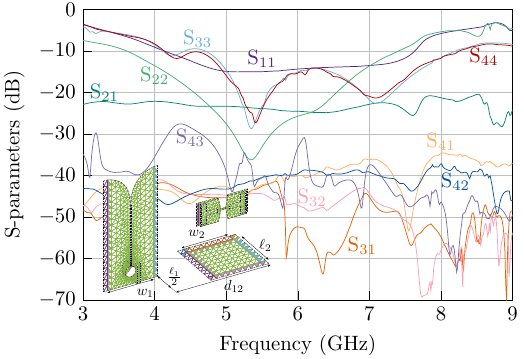}
  \end{center}
  \caption{Selected S-parameters of the two closely placed $9\times 9$ array
    antennas, calculated using the proposed array decomposition method with the
    displayed nine-component representation, where $w_1 = w_2 =$
    \SI{20.5}{\milli\meter}, $\ell_1 = \ell_2 =$ \SI{20.5}{\milli\meter}, and
    $d_{12} =$ \SI{225.5}{\milli\meter}.}
  \label{fig:Sparams}
\end{figure}

In this example, the scattering matrix is calculated for two closely placed
$9\times 9$ array antennas (depicted in Fig.~\ref{fig:2arraysGeo}) in the
frequency band \SI{3}{\giga\hertz} to \SI{9}{\giga\hertz}. The first array
consists of T-slot loaded dipole elements \cite{kolitsidas_rectangular_2014},
whereas the second array consists of Vivaldi elements. By constructing array
components according to the inset in Fig.~\ref{fig:Sparams}, a representation of
the array depicted in Fig.~\ref{fig:2arraysGeo} is obtained. Note that the
Vivaldi unit cell is only electrically connected in the $y$-direction, whereas
the T-slot loaded dipole unit cell is connected in both $x$ and $y$-direction
through the ground plane. 

The array representation in the inset of Fig.~\ref{fig:Sparams} consists of two
disjoint array representations, one for each type of array. Upon imposing the
dimension requirement that the element component of the two representations
needs to be of the same size, \ie, $w_1 = w_2$, and $\ell_1 = \ell_2$, the two
disjoint arrays can be considered as "one" array with the nine component
representation.  This methodology enables the calculations of a number of arrays
situated closely, each with a different array element. 

Using the proposed array representation, \SI{9.53}{\giga\byte} is necessary to
represent the \gls{MoM} impedance matrix. If instead the arrays in
Fig.~\ref{fig:2arraysGeo} are meshed directly, \num{157 220} RWG edge elements are
obtained, resulting in a \gls{MoM} impedance matrix that requires
\SI{395}{\giga\byte} in memory.

To obtain the mutual coupling in the given frequency band, with a maximum
relative error of $10^{-2}$, \num{36} scattering matrices are calculated at
frequency points determined iteratively using the Theta I algorithm
\cite{akerstedt_adaptive_2025}. Between the \num{36} samples, interpolation
with the Loewner framework \cite{lefteriu_new_2010} is used to approximate the
response. The procedure yields a set of \num{500} scattering matrices, each with
a size of $162\times 162$, corresponding to equidistantly spaced frequency
points in the range \SI{3}{\giga\hertz} to \SI{9}{\giga\hertz}. Because of the
large number of ports, only the scattering parameters corresponding to the
highlighted ports in Fig.~\ref{fig:2arraysGeo} are plotted in
Fig.~\ref{fig:Sparams}.

\subsection{Far-field of a $32\times32$ BoR Array}
\begin{figure}[ht!]
  \begin{center}
    \subfloat[]{\includegraphics[width=62pt]{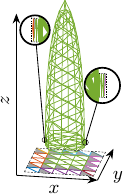}}
    \hfill
    \subfloat[]{\includegraphics[width=184pt]{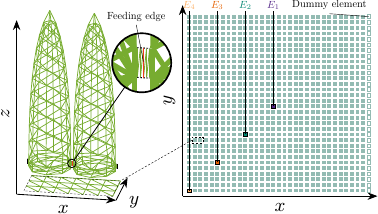}}
  \end{center}
  \caption{BoR array representation. (a) The nine-component representation. (b)
  Schematic view of the array where the connection between two element components
  is displayed. Additionally, array elements whose \gls{EEP} are plotted
  are highlighted.}\label{fig:EEPTesting}
\end{figure}

In this example, the methodology of Section \ref{sec:farfield} is utilized to
calculate the far-field of a $32\times 32$ array consisting of BoR
elements\cite{holter_dual-polarized_2007}. Here, the array components depicted
in Fig.~\ref{fig:EEPTesting}(a) are used to represent the array. In
Fig.~\ref{fig:EEPTesting}(b), the array is visualized, where the
connection between elements is highlighted. The feeding edge is located in the
right connecting strip in the array element component. Subsequently, the
rightmost column of the array consists of dummy elements. With the proposed
array representation \SI{3.97}{\giga\byte} is necessary to represent the
\gls{MoM} impedance matrix. Meshing the array with traditional methods yields
\num{350 519} RWG edge elements, resulting in a \gls{MoM} impedance matrix that
requires \SI{1.97}{\tera\byte} in memory.

\begin{figure}[ht!]
  \begin{center}
    \begin{tikzpicture}
      \node[inner sep=0pt]  at (0pt,0)
    {\includegraphics
    [width=252pt]{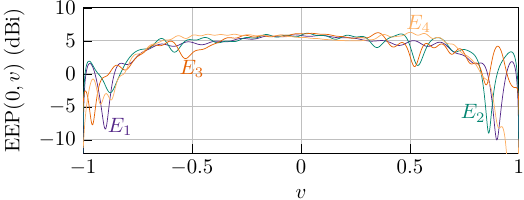}
    };

      \node[inner sep=0pt]  at (0pt,-100pt)
    {\includegraphics
    [width=252pt]{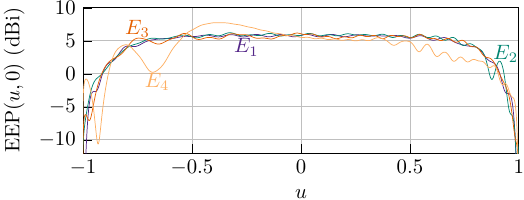}
    };

    \end{tikzpicture}
  \end{center}
  \caption{Co-polarized \gls{EEP} of the four elements highlighted in
  Fig.~\ref{fig:EEPTesting}, given for the H-plane (upper plot) and the E-plane
(lower plot).}\label{fig:EEPEH}
\end{figure}

The \gls{EEP} of the four elements highlighted in Fig.~\ref{fig:EEPTesting}(b)
are calculated and plotted in Fig.~\ref{fig:EEPEH}. Here, we use Ludwig's third
definition of co- and cross-polarization \cite{ludwig_definition_1973}.

Observing the \gls{EEP}s in both the E- and H-plane, the interior array elements
yield similar patterns, whereas the pattern of the corner array element differs
noticeably, as expected. 

\begin{figure}[ht!]
  \begin{center}
    \begin{tikzpicture}
      \node[inner sep=0pt]  at (0pt,0)
    {\includegraphics
    [width=252pt]{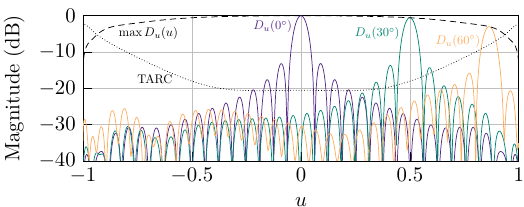}
    };

      \node[inner sep=0pt]  at (0pt,-100pt)
    {\includegraphics
    [width=252pt]{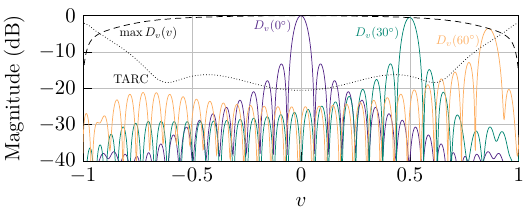}
    };

    \end{tikzpicture}
  \end{center}
  \caption{Co-polarized far-field of the $32\times 32$ BoR array, given in the
    H-plane (upper plot) and in the E-plane (lower plot).
}\label{fig:FFtot}
\end{figure}

Calculating all the \gls{EEP}s of the $32\times32$ BoR array in addition to an
appropriately chosen excitation vector $a$, the directivity $D(u,v)$ is
obtained. In Fig.~\ref{fig:FFtot}, the normalized directivity $D_v(v)$ and
$D_u(u)$ are plotted in the E- and H-plane, respectively. Furthermore,
calculating all the \gls{EEP}s yields all the port currents from which the
scattering matrix $S$ can be calculated. With the scattering matrix $S$ and the
excitation vector $a$, the \gls{TARC} $\Gamma_a^t$ is calculated as
\cite{manteghi_multiport_2005}:
\begin{equation}
  \begin{aligned}
    \Gamma_a^t = \left(\frac{\sum_{i = 1}^{p}|b_i|^2}{\sum_{i = 1}^{p}|a_i|^2}
    \right)^{\frac{1}{2}} \text{,} \quad b = Sa \text{.}
\end{aligned}
  \label{eq:TARC}
\end{equation}
In Fig.~\ref{fig:FFtot}, the \gls{TARC} is plotted and observed to remain below
$-10$ dB in the range $\theta = \pm 52.9^{\circ}$ for the E-plane, and $\theta =
\pm 52.5^{\circ}$ for the H-plane. Additionally, the broadside directivity of
the $32\times32$ BoR array reaches a level of $35$ dBi.

\subsection{Scaling of the Algorithm}
\begin{figure}[ht!]
  \begin{center}
    \includegraphics[width=\columnwidth]{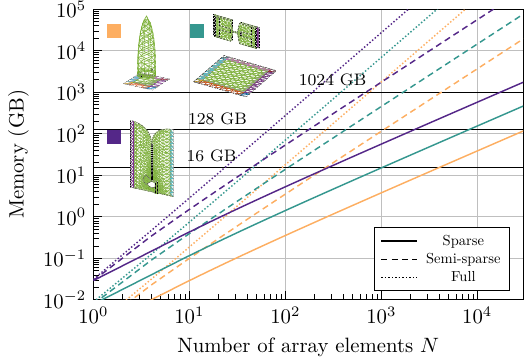}
  \end{center}
  \caption{Memory to allocate for the showcased $N= N_x\times N_y$ arrays, for
    the matrix representations \emph{Full}, \emph{Semi sparse}, and
    \emph{Sparse}.
}\label{fig:scalingAlgo}
\end{figure}

The proposed array decomposition method reduces the allocated memory by an order
of magnitude. Fig.~\ref{fig:scalingAlgo} depicts the memory allocation for an $N
= N_x \times N_y$ array, for the three showcased array representations. Here,
the three matrix representations, \emph{Full}, \emph{semi-sparse}, and
\emph{sparse}, are defined as:
\begin{itemize}
  \item \emph{Full}: The complete matrix $Z^{\mathrm{part}}$,
    $\mathcal{O}(N_x^2N_y^2)$
  \item \emph{Semi-sparse}: The first level of block Toeplitz structure of
    matrix $A$ and complete matrices $B$ and $C$, \ie,
    Fig.~\ref{fig:Amatrices}(b), ~\ref{fig:Bmatrices}(a),
    and~\ref{fig:Cmatrices}(a), $\mathcal{O}(N_x^2N_y), N_y>N_x$
  \item \emph{Sparse}: The second level of block Toeplitz structure of $A$ and
    the sparse representations of $B$ and $C$, \ie, Fig.~\ref{fig:Amatrices}(c),
    ~\ref{fig:Bmatrices}(b), and~\ref{fig:Cmatrices}(b), $\mathcal{O}(N_xN_y)$
\end{itemize}

While the \gls{MoM} impedance matrix can be represented with only a few
submatrices (sparse representation), some inversion algorithms may require the
first level of block Toeplitz structure (semi-sparse representation), \eg, the
matrix-valued Rybicki algorithm.

\section{Conclusions}
\label{sec:Conclusions}
In this article, we have proposed an array decomposition method for a
memory-efficient calculation of large electrically \emph{connected} array
antennas, using RWG \gls{MoM}, that scales as $\mathcal{O}(N_xN_y)$ for an
$N_x\times N_y$ array [as opposed to $\mathcal{O}(N_x^2N_y^2)$]. The method
relies on a nine-component representation of a finite array antenna (eight
margin components and a center array element component), such that the majority
of the \gls{MoM} impedance matrix is obtained as a multilevel block Toeplitz
matrix, enabling a great reduction in necessary submatrix calculations. Here, we
have explicitly shown the necessary submatrix calculations.

The proposed method can be used for the calculation of \gls{EEP}, the full
scattering matrix, the surface currents, and the far-field of a large
electrically connected array. Additionally, the calculation of the far-field is
greatly accelerated using the nine-component representation, as the spatial
Fourier transform integral is only computed once on the nine components'
apertures.

Lastly, the nine-component representation can be used to represent multiple
disjoint and displaced arrays with different array elements, as demonstrated.
Cases such as the mutual coupling between two array antennas (each with a
different array element) closely placed can be calculated rapidly using this
strategy.

\appendices
\section{Memory Allocation for the Nine Components Representation}
Let $e_i$, $i\in[1,9]$ denote the number of edge elements in array component
$i$. The \gls{MoM} impedance matrix of an $N_x \times N_y$ array then requires
an allocation of $MEM(N_x, N_y) = MEM_A + MEM_B + MEM_C$ memory units, where 
\begin{equation}
  \begin{aligned}
    MEM_A(N_x,N_y) = \left(2N_x N_y - N_x - N_y + 1\right)e_5^2 \text{,}
  \end{aligned}
  \label{eq:calcScalingA}
\end{equation}
\begin{equation}
  \begin{aligned}
    &MEM_B(N_x,N_y) = (2N_xN_y-N_x)e_5(e_2 + e_8)\\
    + &(2N_xN_y-N_y)e_5(e_6 +  e_4)+ N_xN_ye_5(e_3+e_9+e_1+e_7) \text{,}
  \end{aligned}
  \label{eq:calcScalingB}
\end{equation}
and
\begin{equation}
  \begin{aligned}
    &MEM_C(N_x,N_y) = N_xN_y (e_2 + e_8 )( e_4 + e_6) + N_y(e_2^2 + e_8^2)\\
&+ N_x(e_6^2 + e_4^2)+ (2N_y-1)e_2e_8 + (2N_x-1)e_6e_4  \\
&+\big(N_x(e_6+e_4)+ N_y(e_2+e_8) \big)(e_3+e_9+e_1+e_7) +e_7^2\\
&+ e_3(e_3 + e_9 + e_1 + e_7) + e_9(e_9 + e_1 + e_7)+ e_1(e_1 + e_7)  \text{.}
  \end{aligned}
  \label{eq:calcScalingC}
\end{equation}


\section*{Acknowledgment}
This work is supported by project nr ID20-0004 from the Swedish Foundation for
Strategic Research, and the Swedish Research Council's Research Environment
grant (SEE-6GIA 2024-06482) for research on sixth-generation wireless systems
(6G), which we gratefully acknowledge.

\ifCLASSOPTIONcaptionsoff
  \newpage
\fi

\bibliographystyle{IEEEtran}
\bibliography{./bibtex/main}

\end{document}

%% file: commands.tex
\usepackage{Carrickc,lettrine}

\makeatletter
\newcommand{\DeclareLatinAbbrev}[2]{%
  \DeclareRobustCommand{#1}{%
    \@ifnextchar{.}{\textit{#2}}{%
      \@ifnextchar{,}{\textit{#2.}}{%
        \@ifnextchar{!}{\textit{#2.}}{%
          \@ifnextchar{?}{\textit{#2.}}{%
            \@ifnextchar{)}{\textit{#2.}}{%
              {\textit{#2.,\ }}}}}}}}%
}
\makeatother
\DeclareLatinAbbrev{\Eg}{E.g}
\DeclareLatinAbbrev{\Ie}{I.e}
\DeclareLatinAbbrev{\etc}{etc}
\DeclareLatinAbbrev{\etal}{et~al}

\newcommand{\ie}{i.e.}
\newcommand{\eg}{e.g.}

\newacronym{MoM}{MoM}{method of moments}
\newacronym{EFIE}{EFIE}{electric-field integral equation}
\newacronym{MLFMM}{MLFMM}{multilevel fast multipole method}
\newacronym{ACA}{ACA}{adaptive cross approximation}
\newacronym{FMM}{FMM}{fast multipole method}
\newacronym{HO-ADM}{HO-ADM}{higher-order array decomposition method}
\newacronym{SFX}{SFX}{synthetic-functions approach}
\newacronym{CBFM}{CBFM}{characteristic basis function method}
\newacronym{SED}{SED}{sub-entire-domain}
\newacronym{DEMCEM}{DEMCEM}{direct evaluation method in computational
electromagnetics}
\newacronym{GMRES}{GMRES}{generalized minimal residual method}
\newacronym{EEP}{EEP}{embedded element pattern}
\newacronym{TARC}{TARC}{total active reflection coefficient}
\newacronym{DDM}{DDM}{domain decomposition method}

%% file: MATLAB_colors.tex
\definecolor{White}{RGB}{255,255,255}

\definecolor{Black}{RGB}{0, 0, 0}

\definecolor{Mblue}{HTML}{0072BD}

\definecolor{Mred}{HTML}{D95319}

\definecolor{Morange}{HTML}{EDB120}

\definecolor{Mpurple}{HTML}{7E2F8E}

\definecolor{Mgreen}{HTML}{77AC30}

\definecolor{Mlightblue}{HTML}{4DBEEE}

\definecolor{Mwine}{HTML}{A2142F}

\definecolor{Corange}{HTML}{E66101}

\definecolor{Cpurple}{HTML}{5E3C99}

\definecolor{meshYellow}{HTML}{FCE469}
\definecolor{meshYellowSep}{HTML}{C8A704}

\definecolor{meshGreen}{HTML}{164F44}

\definecolor{portGreen}{HTML}{004d00}

\definecolor{makarovYellow}{HTML}{FFBF3B}

\definecolor{makarovYellow2}{HTML}{E8D792}

\definecolor{makarovBlue}{HTML}{93BBB3}

\definecolor{makarovGreen}{HTML}{014D43}

\definecolor{makarovDarkGreen}{HTML}{073727}

\definecolor{newOrange}{HTML}{fec44f}

\definecolor{Mmagenta}{HTML}{cc64be}
\definecolor{Mlightgreen}{HTML}{79fc80}
\definecolor{Mteal}{HTML}{5a8d8f}

\definecolor{portColor}{HTML}{f53c5e}